\documentclass[12pt]{amsart}
\usepackage{amscd,amssymb}
\usepackage[graph,frame,poly,arc]{xy}  
\usepackage[plainpages,backref,urlcolor=blue]{hyperref}
\usepackage{mathtools}

\theoremstyle{plain}
\newtheorem{thm}[subsection]{Theorem}
\newtheorem{lem}[subsection]{Lemma}

\newtheorem{cor}[subsection]{Corollary}

\theoremstyle{definition}

\newtheorem{ex}[subsection]{Example}

\numberwithin{equation}{section}
\newcommand{\M}{{\mathcal M}}

\newcommand{\al}{{\alpha}}

\newcommand{\C}{\mathbb{C}}
\newcommand{\PP}{\mathbb{P}}

\newcommand{\dd}{{\rm d}}

\begin{document}

\title [On the Jacobian syzygies for generic toric models]
{On the Jacobian syzygies for generic toric models}

\author[Alexandru Dimca]{Alexandru Dimca}
\address{Universit\'e C\^ ote d'Azur, CNRS, LJAD, France and Simion Stoilow Institute of Mathematics,
P.O. Box 1-764, RO-014700 Bucharest, Romania}
\email{Alexandru.Dimca@univ-cotedazur.fr}

\author[Gabriel Sticlaru]{Gabriel Sticlaru}
\address{Faculty of Mathematics and Informatics,
Ovidius University
Bd. Mamaia 124, 900527 Constanta,
Romania}
\email{gabriel.sticlaru@gmail.com }



\begin{abstract} 
To a generic hypersurface in the affine torus $(\C^*)^n$ we associate a hypersurface arrangement in the projective space $\PP^n$ consisting of the $n+1$ coordinate hyperplanes and a generic hypersurface, and compute the minimal graded resolutions of the corresponding Jacobian algebra.
\end{abstract}
 
\maketitle


\section{Introduction} 

Let $S=\C[x_0,x_1, \ldots, x_n]$ be the polynomial ring in $n+1 \geq 3$ variables $x_0,x_1, \ldots, x_n$ with complex coefficients, and let $V:f=0$ be a reduced hypersurface of degree $d\geq 3$ in the complex projective plane $\PP^n$. 
We denote by $J_f$ the Jacobian ideal of $f$, i.e., the homogeneous ideal in $S$ spanned by the partial derivatives $f_j$ of $f$ with respect to $x_j$ for $j=0, \ldots,n$, and  by $M(f)=S/J_f$ the corresponding graded quotient ring, called the Jacobian (or Milnor) algebra of $f$.
Consider the graded $S$-module of Jacobian syzygies of $f$ or, equivalently, the module of derivations killing $f$, namely
\begin{equation}
\label{eqD0}
D_0(f)= \{\theta \in Der(S)  :  \theta(f)=0\}= \{\rho=(a_0,a_1, \ldots, a_n) \in S^{n+1}  :   \sum_{j=0}^na_jf_j=0\}.
\end{equation}
We say that $V:f=0$ is an {\it $m$-syzygy hypersurface} if  the module $D_0(f)$ is minimally generated by $m$ homogeneous syzygies, say $\rho_1,\rho_2,\ldots ,\rho_m$, of degrees $d_j=\deg \rho_j$ ordered such that $$d_1\leq d_2 \leq \ldots \leq d_m.$$ 
We call these degrees $(d_1, \ldots, d_m)$ the {\it exponents} of the hypersurface $V$.

There is an increasing interest in finding the minimal graded resolution of the Jacobian algebra $M(f)$, or at least the exponents of the graded $S$-module $D_0(f)$, see for instance \cite{Bernd1, Bernd2}, where connections with {\it  likelihhood correspondence and applications to statistics and physics are discussed}. Note for instance that the graded $S$-module $D_0(f)$ is denoted by $S_0(V)$ in \cite[Lemma 2.2 (2)]{Bernd1}.

The simplest case is when $V:f=0$ is a smooth hypersurface in $\PP^n$, since then the partial derivatives $f_0, \ldots, f_n$ form a regular sequence in $S$ and hence a minimal resolution for the corresponding Jacobian algebra $M(f)$ is given by 
\begin{equation}
\label{res1}
 0 \to S(-e_{n+1})^{c_{n+1}} \xrightarrow{u_n} S(-e_{n}) ^{c_{n}}\xrightarrow{u_{n-1}} \ldots  \xrightarrow{u_1} S(-e_{1})^{c_{1}}\xrightarrow{u_0} S(-e_{0})^{c_{0}},
\end{equation}
where $c_k=\binom{n+1}{k}$ and $e_k=k(d-1)$ for $k=0, \ldots,n+1$.
Moreover, the morphisms 
$$u_k:S(-e_{k+1})^{c_{k+1}} \to S(-e_{k})^{c_{k}}$$ correspond under an obvious identification
$$S^{c_{k+1}}=\Omega^{n-k} (k-n)\text{ and } S^{c_{k}}=\Omega^{n+1-k} (k-n-1)$$
to the cup product by the differential $\dd f$ where $\Omega^s$ denotes the differential s-forms on $\C^{n+1}$ with polynomial coefficients, see
\cite[Equation (6.2.5) and Example (6.2.13)]{STH}.
In particular, a smooth hypersurface of degree $d$ in $\PP^n$ is an
$N$-syzygy hypersurface with exponents
$$d_1= \ldots =d_N=d-1 \text{ where } N=\binom{n+1}{2}.$$
As soon as singularities arise, the problem becomes much more complicated, as one can see by looking at the extensive literature for curves, or even line arrangements, in $\PP^2$, see for instance \cite{Mich,expo,minTjurina,Sch0,To}.
One case where the answer is known is when $V:f=0$ is a normal crossing arrangement of $d>n+1$ hyperplanes in $\PP^n$, see \cite[Corollary 4.5.4]{RT} or \cite{Y}. In this case the minimal resolution for $M(f)$ has the same shape as in \eqref{res1}, but this time
\begin{equation}
\label{res2A}
c_0=1, \ c_1=n+1, \ c_k=\binom{d+k-n-4}{k-2}\binom{d-1}{n+1-k} 
\end{equation}
and 
\begin{equation}
\label{res2B}
e_0=0, \ e_1=d-1, \ e_k=2d+k-n-3
\end{equation}
for $2 \leq k \leq n+1$.
When $V$ has all the irreducible components smooth of degree $>1$ and they satisfy a rather involved genericity condition, the corresponding minimal resolution is described in \cite[Theorem 3.6]{BRSS}.

In this note we consider the following geometric setting. Let $H_j: \ell_j=0$ for $j=0,\ldots,n$ be $n+1$ hyperplane in $\PP^n$ such that
\begin{equation}
\label{CH}
H_0 \cap H_1 \cap \ldots \cap H_{n+1}=0,
\end{equation}
and let $W:g=0$ be a smooth degree $e$ hypersurface in $\PP^n$.
With this notation, we have the following result.
\begin{thm}
\label{thm1}
Let $V=W \cup H_0 \cup H_1 \cup \ldots \cup H_n$ be the projective hypersurface in $\PP^n$ defined by $f=g\ell_0 \ell_1 \ldots \ell_n=0$.
We assume that $V$ is a normal crossing divisor.
Then the minimal resolution of the Jacobian algebra $M(f)$ has the shape given in \eqref{res1} with
$$c_k=\binom{n+1}{k}, \ e_0=0, \ e_1=e+n, \text{ and } e_k=ke+n+1$$ for $2 \leq k \leq n+1$.
In particular,  $V$ is an
$N$-syzygy hypersurface with exponents
$$d_1= \ldots =d_N=e+1 \text{ where } N=\binom{n+1}{2}.$$
Moreover, for a fixed collection of hyperplanes $H_0, \ldots, H_n$ as above, there is a Zariski open subset of polynomials $g$ in $S_e$ for which the above assumption holds.
\end{thm}
Note that up-to a change of coordinates we can suppose that $\ell_j=x_j$ for all $j=0, \ldots, n$ and hence the complement
$$\M(V)=\PP^n \setminus V$$
coincides with the complement
$$\M(W')=(\C^*)^n \setminus W',$$
where we identify $(\C^*)^n$ to $\PP^n \setminus (H_0 \cup \ldots \cup H_n)$ by setting $x_0=1$ and the hypersurface $W'$ is defined in 
$(\C^*)^n$ by the equation
$$g'(x_1,x_2,\ldots,x_n)=g(1,x_1,x_2, \ldots,x_n).$$
To associated  a graded  object to the Jacobian ideal of such a hypersurface $W'$, passing to the Jacobian algebra $M(f)$ described above seems to be the natural approach, see also \cite[Formula (8)]{Bernd1} and the discussion around it concerning a generic torus model. 
For background on likelihood geometry and toric models we refer the interested reader to \cite{BCF}. This view-point explains the title of this note.

The following result is a consequence of the proof of Theorem \ref{thm1}.
\begin{cor}
\label{cor1}
Let $W: g=0$ be a smooth projective hypersurface in $\PP^n$ of degree $e>0$, and set $H_i: x_i=0$ for any $i=0, \ldots, x_n$.
We assume that 
$$V:f=x_0\ldots x_n g=0$$
 is a normal crossing divisor.
Set  $g_i'=x_ig_i+g$ for $i=0, \ldots, x_n$, where $g_i$ is the partial derivative of $g$ with respect to $x_i$.
Then a minimal set of generators for the graded $S$-module $D_0(f)$ is given by the syzygies
$$\rho_{ij}' \in S_{e+1}^{n+1}$$
for $0 \leq i <j \leq n$, having only two non-zero components, namely the $i$-th component equal to $x_ig'_j$ and the $j$-th component equal to $-x_jg'_i$.
\end{cor}
For more general results of this type the reader is referred to 
\cite[Theorem 2.11]{Bernd1} and \cite[Proposition 5.1 and Corollary 5.2]{Bernd2}. Note that the isomorphisms provided by these results are fairly indirect, and so not easy to use in practice.

In the next section we prove Theorem \ref{thm1} and we show  that the normal crossing assumption in this Theorem is necessary, see Example \ref{ex1}. 

\medskip

We would like to thank Bernd Sturmfels for initial motivation and very useful exchanges concerning this note.

\section{Proof of Theorem \ref{thm1}} 
We may and do choose the coordinates on $\PP^n$ such that
$\ell_j=x_j$ for all $j=0, \ldots, n$. In this new coordinates, we have
$$f=gh, \text{ where } h=x_0 x_1 \ldots x_n.$$
It follows that 
$$f_i=g_ih+gh_i=h_i(x_ig_i+g).$$
Suppose we have a syzygy
\begin{equation}
\label{syz1}
a_0f_0+a_1f_1+ \ldots + a_nf_n=0.
\end{equation}
For any fixed index $i$, since $g$ is not divisible by $x_i$ and all of the $h_j$ for $j \ne i$ are divisible by $x_i$, it follows that $a_i$ is divisible by $x_i$. Therefore we have $a_i=x_iA_i$, for some polynomials $A_i \in S$. If we replace this equality in \eqref{syz1} and divide by the common factor $h=x_ih_i$, we get the following equality
\begin{equation}
\label{syz2}
A_0(x_0g_0+g)+A_1(x_1g_1+g)+ \ldots + A_n(x_ng_n+g)=0.
\end{equation}
Assume that $g'_i=x_ig_i+g$ form a regular sequence in $S$. Then the module of the corresponding syzygies
$$D'=\{(A_0, \ldots,A_n) \ : \ A_0g'_0+\ldots +A_ng'_n=0\}$$
is generated by the Koszul type syzygies
$$\rho_{ij}\in S_e^{n+1}$$
for $0 \leq i <j \leq n$, having only two non-zero components, namely the $i$-th one equal to $g'_j$ and the $j$-th one equal to $-g'_i$.
It follows that the $S$-module $D_0(f)$ is generated by the syzygies
$$\rho_{ij}' \in S_{e+1}^{n+1}$$
for $0 \leq i <j \leq n$, having only two non-zero components, namely the $i$-th one equal to $x_ig'_j$ and the $j$-th one equal to $-x_jg'_i$.
This shows that the minimal resolution of $M(f)$ starts as follows
$$  S(-(2e+n+1))^{c_2} \to S(-(e+n))^{c_1} \to S.$$
The image of the first morphism is clearly just $D_0(f)(-(e+n)))$.
On the other hand, the morphism
$$\al: D'(-1) \to D_0(f), \  \rho_{ij} \mapsto \rho_{ij}',$$
induces an isomorphism of graded $S$-modules $D'(-1) \to D_0(f)$.
Since the $g_i'$ are supposed to form a regular sequence, the quotient graded $S$-module $S/(g_0',\ldots,g_n')$ has a minimal resolution of the shape given in \eqref{res1} with
 $c_k=\binom{n+1}{k}$ and $e_k''$ which replaces $e_k$ given by $e_k''=ke$ for $k=0, \ldots,n+1$, see for instance \cite[Corollary 17.5]{Eis}.
 It follows that the graded $S$-module $D'$, which is the kernel of the morphism
 $$S^{c_{1}} \xrightarrow{u_0} S$$
in this minimal resolution, has a minimal resolution given by
 $$0 \to S(e_1''-e_{n+1}'')^{c_{n+1}} \xrightarrow{u_n} S(e_1''-e_{n}'') ^{c_{n}}\xrightarrow{u_{n-1}} \ldots  \xrightarrow{u_1} D' \to 0.$$ 
Twisting this exact sequence by $(-1)-(e+n)=-(e+n+1)$ and noticing that
$$e_1''-e_k''-(e+n+1)=e-ke-e-n-1=-(ke+n+1),$$
we get the value of $e_k$ for $2 \leq k \leq n+1$ claimed in Theorem \ref{thm1}.

It remains to show the following result.
\begin{lem}
\label{lem1}
Let $W: g=0$ be a smooth projective hypersurface in $\PP^n$ of degree $e>0$, and set $H_i: x_i=0$ for any $i=0, \ldots, x_n$.
We assume that $V:f=x_0\ldots x_n g=0$ is a normal crossing divisor.
Then the polynomials $g_i'=x_ig_i+g$ for $i=0, \ldots, x_n$ form a regular sequence in $S$, where $g_i$ is the partial derivative of $g$ with respect to $x_i$.
Moreover, there is a Zariski open subset of polynomials $g$ in $S_e$ for which the above assumption holds.
\end{lem}
\proof

Let $W:g=0$ be a smooth hypersurface of degree $e>0$ in $\PP^n$. For each proper subset $I$ in $\{0,1,\ldots,n\}$, we consider the corresponding edge
$$E_I=\cap_{i \in I}H_i$$
and we note the $\dim E_I=n-|I|.$ 
Moreover, when $\dim E_I >0$, the intersection $W_I=W\cap E_I$ is a smooth hypersurface in $E_I$ if and only if $E_I$ is not tangent to $W$ at some point $p \in W_I$, in other words $E_I$ is not contained in the projective tangent space $T_pW$ to $W$ at $p$. By convention, if $\dim E_I=0$ and the point $E_I$ is in $W$, we say that $E_I$ is tangent to $W$, since $E_I \subset T_pW$ in such a case.
It is clear that $V$ is a normal crossing divisor if and only if every edge $E_I$ is not tangent to $W$.

The set of $n+1$ homogeneous polynomials $g_i'$ of degree $e>0$ form a regular sequence in $S$ if and only if the system of equations
$$g_0'(x)=g_1'(x)= \ldots = g_n'(x)=0$$
admits only the trivial solution in $\C^{n+1}$, see for instance \cite[Proposition 7.23]{RCS}. Assume that there is a solution $p \ne 0$ for this system. Then by adding these $n+1$ equations and using the Euler relation, we get $g(p)=0$, that is the corresponding point 
$$p=(p_0:p_1: \ldots :p_n) \in \PP^n$$ is in fact in $W$. We may assume, by reordering the variables $x_j$'s if necessary, that there is an integer $k\geq 0$ such that $p_0,p_1, \ldots,p_k$ are all non-zero, and 
$$p_{k+1}=\ldots = p_n=0.$$
It follows that $p \in E_I$, where $I=\{k+1,\ldots,n\}$.
We have seen that $g(p)=0$, and hence the above system implies that
$$g_0(p)=\ldots=g_k(p)=0.$$
Hence the tangent vector space $T_pW$ to $W$ at the point $p$ is given by
$$g_{k+1}(p)x_{k+1}+ \ldots +g_n(p)x_n=0.$$
Since the tangent space $T_pE_I=E_I$ to $E_I$ at the point $p$ is given by
$$x_{k+1}=\ldots =x_n=0,$$
it follows that $T_pE_I \subset T_pW$, that is $E_I$ is tangent to $W$ at the point $p$. This contradiction proves our first claim in Lemma \ref{lem1}. 

To prove the last claim there, note that the assumption in 
Lemma \ref{lem1} holds for the Fermat type polynomial
$$g=x_0^e+x_1^e+ \ldots +x_n^e.$$
On the other hand, this condition is an open condition in the Zariski topology of $S_e$. The most rapid way to see this is to use the fact that $V$ is a normal crossing divisor if and only if the Euler discriminant defined in \cite[Section 4]{Bernd2} is nonzero at $V$.
An alternative proof for this result can be obtained using a transversality approach, see for instance \cite{Kl} or \cite[Section 2 in Chapter 2]{Gi} and \cite[Lemma 2.4.6 and Remark on p. 54]{GG}.
\endproof


\begin{ex}
\label{ex1} Consider one of the simplest cases of Theorem \ref{thm1}, namely when $n=e=2$, for instance when
$$V:f=x_0x_1x_2(x_0^2+x_1^2+x_2^2)=0,$$
namely $V$ is the union of a triangle and a smooth conic transversal to the triangle in 6 distinct points.
Then the minimal resolution of $M(f)$ given by Theorem \ref{thm1} has the following shape
$$0 \to S(-9) \to S^3(-7) \to S^3(-4) \to S.$$
In \cite[Example 4.14]{DP} it is shown that if $V':f'=0$ is a union of a triangle and a smooth conic which is either tangent to the 3 sides of the triangle, or passes through the 3 vertices of the triangle, then the minimal resolution of the corresponding Jacobian algebra $M(f')$ has the shape
$$0 \to S^2(-6) \to S^3(-4) \to S.$$
Similarly, the condition \eqref{CH} on the hyperplanes $H_i$ is necessary, that is the defining equations should be linearly independent.
The curve
$$C:f=(x_0-x_1)(x_1-x_2)(x_0-x_2)(x_0^2+x_1^2+x_2^2)=0,$$
where this condition is not fulfilled, has as minimal resolution for $M(f)$
the following sequence
$$0 \to S(-9) \to S(-6)\oplus S(-7) \oplus S(-8) \to S(-4)^3 \to S,$$
obtained using SINGULAR, and which is quite different from the resolution given in Theorem \ref{thm1}.
Note that in this case the 3 lines are transversal to the conic $x_0^2+x_1^2+x_2^2=0$ and there are only the triple point $(1:1:1)$ where the normal crossing condition fails. This curve $C$ has as exponents $(d_1,d_2,d_3)=(2,3,4)$, and hence it is a minimal plus-one generated curve as defined in \cite{Brian}. In fact, one can obtain the exponents of this curve $C$ without using SINGULAR by first showing that $d_1=2$ and then applying \cite[Theorem 1.5]{Brian}.
It follows that some assumption is necessary in Theorem \ref{thm1} in order to have a minimal resolution for $M(f)$ of the stated shape.
\end{ex}

\end{document}